\documentclass[11pt]{amsart}

\usepackage{amsthm}
\usepackage{amssymb}
\usepackage{graphicx}									 
\usepackage{stmaryrd}                  
\usepackage{footmisc}                  
\usepackage{paralist}
\usepackage{wrapfig}
\usepackage[draft]{pdfcomment}
\usepackage{bbm}
\usepackage[T1]{fontenc}
\usepackage[utf8]{inputenc}
\usepackage[arrow, matrix, curve]{xy}
\usepackage{color}

\usepackage{tikz-cd}

\usepackage{enumitem}

\newtheorem{thm}{Theorem}[section]

\newtheorem{prp}[thm]{Proposition}

\theoremstyle{definition}

\newtheorem{qst}[thm]{Question}

\theoremstyle{remark}

\newcommand{\Z}{\mathbb{Z}}

\newcommand{\Orb}{\mathcal{O}}

\title[Good, but not very good orbifolds]{Good, but not very good orbifolds}

\author[C. Lange]{Christian Lange}
\address{Ludwig-Maximilians-Universit\"at M\"unchen, Mathematisches Institut\newline\indent Theresienst. 39, 80333 München, Germany}
\email{lange@math.lmu.de}


\begin{document}

\begin{abstract}
We construct examples of (effective) closed orbifolds which are covered by manifolds, but not finitely so.
\end{abstract}

\maketitle

Orbifolds were introduced by Satake under the name of V-manifolds \cite{MR0079769} and later rediscovered by Thurston, who generalized the theory of coverings and fundamental groups to the setting of orbifolds \cite{Thurston}. Following Thurston an orbifold is called (\emph{very}) \emph{good} if it is (finitely) covered by a manifold. Otherwise it is called \emph{bad}. A good orbifold can be thought of as a proper action of a discrete group on a manifold \cite{Da11}. The corresponding orbifold is very good if and only if the acting group has a subgroup of finite index such that the action restricted to this subgroup is free.

Up to dimension $3$ any good compact orbifold (admits a smooth structure and) is in fact very good as a consequence of the geometrization theorem for good compact orbifolds, see \cite{BBP05} or \cite{KL14}, and \cite{Sc83} for the elementary $2$-dimensional case. 

In \cite{Fl24} the question is raised, whether any good compact orbifold is very good. Good noncompact orbifolds which are not very good exist in all dimensions $n\geq 2$ and can be easily constructed by including infinitely many local groups with unbounded torsion. There exist nonlinear Lie groups with cocompact lattices that do not have torsion free subgroups of finite index \cite{Ra84}, but these lattices do not act effectively on the corresponding nonpositively curved symmetric spaces, cf. Proposition \ref{prp:loc_sym_good}.

In this note we construct an example of a good compact (effective) smooth orbifold (all of whose local groups are cyclic) which is not very good in any dimension $n\geq 4$.

\section{The example}\label{sec:example}

We start with a finitely presented simple group $\Gamma$ with torsion, for instance, Thompson's group $V$ \cite{CFP96,Hi74}. By Poincaré's theorem the group $\Gamma$ does not have finite index subgroups: any such subgroup $\Gamma'$ would give rise to a proper finite index normal subgroup, namely the subgroup that acts trivially on the cosets $\Gamma/\Gamma'$. For any $n\geq 4$ we can construct a closed $n$-manifold $M$ with fundamental group $\Gamma$ (via surgeries starting with a connected sum of copies of $S^{n-1}\times S^1$), see \cite[Section~52]{ST80} or \cite[Theorem~5.1.1]{CZ93}. Inside $M$ we choose an embedded loop that represents a torsion element $g$ in the fundamental group, say of order $p\geq 2$. Now we cut out an embedded tubular neighborhood of this loop and obtain an $n$-manifold $N$ with boundary $S^{n-2} \times S^1$ and fundamental group $\Gamma$. The inclusion of the boundary into $N$ induces a map on the level of fundamental groups which sends a generator of the fundamental group of $S^1$ to $g$. The group $\Z_p$ acts on the disk $D^2$ by rotations. The $n$-orbifold $S^{n-2}  \times D^2/\Z_p$ has boundary $S^{n-2} \times S^1$ as well and (orbifold) fundamental group $\Z_p$. The inclusion of the boundary into $S^{n-2} \times D^2/\Z_p$ induces a map on the level of fundamental groups which sends a generator of the fundamental group of $S^1$ to a generator of $\Z_p$. We glue $N$ and $S^{n-2}  \times D^2/\Z_p$ along their boundary. By the Seifert--Van Kampen theorem the fundamental group of the resulting orbifold $\mathcal{O}$ is isomorphic to $\Gamma$ as well. Therefore, the orbifold $\mathcal{O}$ does not have any finite covers and is, in particular, not very good. 

We claim that $\mathcal{O}$ is covered by a manifold. The universal covering $\tilde N$ of $N$ is a manifold with boundary each of whose components is homeomorphic to $S^{n-2}\times S^1$ as $g$ has finite order. The subgroup of the deck transformation group that leaves one of these boundary components invariant is isomorphic to $\Z_p$. The universal covering of $S^{n-2}  \times D^2/\Z_p$ is $S^{n-2}  \times D^2$ with boundary $S^{n-2}\times S^1$. The restrictions of these two coverings to a boundary component coincide up to an $\Z_p$-equivariant identification induced by a lift of the gluing map downstairs. Gluing a copy of $S^{n-2}  \times D^2$ to each boundary component of $\tilde N$ via these $\Z_p$-equivariant identifications yields a desired manifold covering of $\Orb$.

\section{Geometric structures and constraints}\label{sec:geometric_constraints}

The work of Cheeger and Gromoll \cite{CG72} implies that a good Riemannian orbifold with nonnegative Ricci curvature is very good, see \cite[Theorem~2.1]{Wi00}. On the other hand, the examples constructed in Section \ref{sec:example} can be equipped with Riemannian metrics of positive scalar curvature via \cite{GL80}.

In nonpositive curvature any complete Riemannian orbifold is good \cite[Theorem~2.15]{BH99}. If the orbifold is in addition locally symmetric, we have
\begin{prp}\label{prp:loc_sym_good} A nonpositively curved complete locally symmetric orbifold is very good.
\end{prp}
\begin{proof} The orbifold is a quotient of a nonpositively curved simply connected symmetric space $M$ by a discrete closed subgroup $\Gamma$ of its isometry group. By discreteness of $\Gamma$ it suffices to find a torsion free subgroup of $\Gamma$ of finite index. 

The isometry group $G$ of $M$ has finitely many connected components, because the connected component of its identity $G^0$ acts transivitely and the isotropy group of a point is compact. Hence, we can assume that $\Gamma$ is contained in $G^0$. By Selberg's lemma \cite[p.~330]{Ra06} it suffices to embed $G^0$ into a matrix Lie group over the complex numbers. The symmetric space $M$ splits as a product of a symmetric space of noncompact type and a Euclidean factor. Its isometry group splits accordingly. Therefore, we can assume that $M$ is of noncompact type. In this case the center of $G^0$ is trivial, which follows for instance from the observation that isotropy groups in $G^0$ have unique fixed points in $M$. 
Hence, $G^0$ embeds into a suitable matrix Lie group via the adjoint representation as desired.
\end{proof}

On the other hand, Kapovich constructs complete negatively curved orbifolds which are not very good \cite{Ka22}, providing a negative answer to a question of Margulis as to whether the conclusion of Selberg's lemma, i.e. the existence of a torsion free subgroup of finite index, holds for all finitely generated discrete isometry groups of Hadamard manifolds \cite{Ma75}. 

\begin{qst} Does there exist a compact nonpositively resp. negatively curved orbifold which is not very good?
\end{qst}

The unit tangent bundles of the $4$-orbifolds from Section \ref{sec:example} with respect to some Riemannian metric give examples of good compact contact $7$-orbifolds, cf. \cite[Theorem 1.5.2]{Ge08}, which are not very good, because their fundamental groups coincide with those of the base by the long exact homotopy sequence and because they are not manifolds as the singularities of the base are not isolated.

\begin{qst} Does there exist a good compact contact $5$-orbifold which is not very good?
\end{qst}

Likewise, we may ask

\begin{qst} Does there exist a good compact (almost) complex resp. symplectic resp. Kähler orbifold (of real dimension 4) which is not very good?
\end{qst}

As for complex orbifolds, examples in complex dimension $3$ might be obtained as twistor spaces \cite{AHS78} from the examples in Section \ref{sec:example} via an orbifold version of the main result of \cite{Ta92}.
\newline
\newline
\emph{Acknowledgements.} I thank Claudio Gorodski and Alexander Lytchak for discussions about symmetric spaces.

\end{document}